\documentclass[11pt, a4paper]{article}
\usepackage{latexsym}
\usepackage{indentfirst}
\usepackage{graphicx}
\usepackage{amsmath}
\usepackage{amssymb}

\begin{document}
\title{Limit cycles by FEM for a one - parameter dynamical
system associated to the Luo - Rudy I model}
\author{C\u{a}t\u{a}lin Liviu Bichir $^{1}$, Adelina
Georgescu $^{2}$, Bogdan Amuzescu $^{3}$, \\ Gheorghe Nistor
$^{4}$, Marin Popescu $^{5}$, Maria-Luiza Flonta $^{6}$,
\\ Alexandru Dan
Corlan $^{7}$, Istvan Svab $^{8}$, \\
1 \ Rostirea Maths Research, Regimentul 11 Siret 27, Gala\c{t}i, Romania, \\
2 \ Academy of Romanian Scientists, \\ Splaiul Independen\c{t}ei 54,
Bucharest, Romania, \\
3, 6, 8 \ Faculty of Biology, University of Bucharest, \\
Splaiul Independen\c{t}ei 91-95, Bucharest, Romania, \\
4, 5 \ University of Pite\c{s}ti, Str. T\^{a}rgul din Vale 1,
Pite\c{s}ti, Romania, \\
7 \ Bucharest University Emergency Hospital, \\ Splaiul
Independen\c{t}ei 169, Bucharest, Romania,
\\ 1 \ catalinliviubichir@yahoo.com, 3 \ bogdan@biologie.kappa.ro, \\
4 \ ghe.nistor@yahoo.com, 5 \ popescumarin67@yahoo.com, \\
6 \ flonta@biologie.kappa.ro, 7 \ alexandru@corlan.net, \\
8 \ istvansvab@icbp.ro}

\date{2010 Oct 10}
\maketitle

\begin{abstract}
An one - parameter dynamical system is associated to the
mathematical problem governing the membrane excitability of a
ventricular cardiomyocyte, according to the Luo-Rudy I model.
Limit cycles are described by the solutions of an extended system.
A finite element method time approximation (FEM) is used in order
to formulate the approximate problem. Starting from a Hopf
bifurcation point, approximate limit cycles are obtained, step by
step, using an arc-length-continuation method and Newton's method.
Some numerical results are presented. 
\end{abstract}

\section{Introduction}
\label{BichirCL_sect_I}

The well-known Hodgkin-Huxley model of the squid giant axon
(\cite{BichirCL_Hodgkin_Huxley}) represented a huge leap forward
compared to earlier models of excitable systems built from
abstract sets of equations or from electrical circuits including
non-linear components, e.g.
\cite{BichirCL_Van_der_Pol_van_der_Mark}. The pioneering work of
the group of Denis Noble made the transition from neuronal
excitability models, characterized by Na+ and K+ conductances with
fast gating kinetics, to cardiomyocyte electrophysiology models, a
field expanding steadily for over five decades
(\cite{BichirCL_Noble7}).  Nowadays, complex models accurately
reproducing transmembrane voltage changes as well as ion
concentration dynamics between various subcellular compartments
and buffering systems are incorporated into detailed anatomical
models of the entire heart (\cite{BichirCL_Noble2007}). The
Luo-Rudy I model of isolated guinea pig ventricular cardiomyocyte
(\cite{BichirCL_LuoRudyI}) was developed in the early 1990s
starting from the Beeler-Reuter model
(\cite{BichirCL_BeelerReuter}). It includes more recent
experimental data related to gating and permeation properties of
several types of ion channels, obtained in the late 1980s with the
advent of the patch-clamp technique
(\cite{BichirCL_Neher_Sakmann}). The model comprises only three
time and voltage-dependent ion currents (fast sodium current, slow
inward current, time-dependent potassium current) plus three
background currents (time-independent and plateau potassium
current, background current), their dynamics being described by
Hodgkin-Huxley type equations. This apparent simplicity, compared
to more recent multicompartment models, renders it adequate for
mathematical analysis using methods of linear stability and
bifurcation theory.

Nowadays, there exist numerous software packages for the numerical
study of finite - dimensional dynamical systems, for example
MATCONT, CL$_{-}$MATCONT, CL$_{-}$MATCONTM
(\cite{BichirCL_MATCONT}, \cite{BichirCL_Cl_MatContM}), AUTO
\cite{BichirCL_Doedel}. In \cite{BichirCL_Kuznetsov},
\cite{BichirCL_Doedel}, \cite{BichirCL_MATCONT},
\cite{BichirCL_Cl_MatContM}, the periodic boundary value problems
used to locate limit cycles are approximated using orthogonal
collocation method. Finite differences method is also considered.
In this paper, limit cycles are obtained for the dynamical system
associated to the Luo-Rudy I model by using finite element method
time approximation (FEM).

\section{Luo-Rudy I model}
\label{BichirCL_sect_II}

The mathematical problem governing the membrane excitability of a
ventricular cardiomyocyte, according to the Luo-Rudy I model
(\cite{BichirCL_LuoRudyI}), is a Cauchy problem
\begin{equation}
\label{BichirCL_e01}
   u(0)=u_{0} \, ,
\end{equation}
for the system of first order ordinary differential equations
\begin{equation}
\label{BichirCL_e02}
   \frac{du}{dt}=\mathcal{F}(\eta,u) \, ,
\end{equation}
where $u$ $=$ $(u_{1},\ldots,u_{8})$ $=$ $(V$, $[Ca]_{i}$, $h$,
$j$, $m$, $d$, $f$, $X)$, $\eta$ $=$ $(\eta_{1},\ldots,\eta_{13})$
$=$ $(I_{st}$, $C_{m}$, $g_{Na}$, $g_{si}$, $g_{Kp}$, $g_{b}$,
$[Na]_{0}$, $[Na]_{i}$, $[K]_{0}$, $[K]_{i}$, $PR_{NaK}$, $E_{b}$,
$T)$, $M=\mathbb{R}^{8}$, $\mathcal{F}:\mathbb{R}^{13} \times M
\rightarrow M$, $\mathcal{F} =
(\mathcal{F}_{1},\ldots,\mathcal{F}_{8})$,
\begin{eqnarray*}
   & \ & \mathcal{F}_{1}(\eta,u)=
      -\frac{1}{\eta_{2}}[I_{st}
         +\eta_{3}{u}_{3}{u}_{4}{u}_{5}^{3}({u}_{1}-E_{Na}(\eta_{7},\eta_{8},\eta_{13}))
         \\
   & \ & \qquad +\eta_{4}{u}_{6}{u}_{7}({u}_{1}-c_{1}+c_{2}\ln{u}_{2})
         \\
   & \ & \qquad +g_{K}(\eta_{10})X_{i}({u}_{1})
            ({u}_{1}-E_{K}(\eta_{7},\eta_{8},\eta_{9},\eta_{10},\eta_{11},\eta_{13})){u}_{8}
         \\
   & \ &  \qquad +g_{K1}(\eta_{10})K1_{\infty}(\eta_{9},\eta_{10},\eta_{13},{u}_{1})
             ({u}_{1}-E_{K1}(\eta_{9},\eta_{10},\eta_{13}))
         \label{BichirCL_e03} \\
   & \ & \qquad +\eta_{5}Kp({u}_{1})({u}_{1}-E_{Kp}(\eta_{9},\eta_{10},\eta_{13}))
          +\eta_{6}({u}_{1}-\eta_{12})] \, ,
         \\
   & \ & \mathcal{F}_{2}(\eta,u)=
      -c_{3}\eta_{4}{u}_{6}{u}_{7}({u}_{1}-c_{1}+c_{2}\ln{u}_{2})
      +c_{4}(c_{5}-{u}_{2}) \, ,
         \\
   & \ & \mathcal{F}_{\ell}(\eta,u)=
      \alpha_{\ell}({u}_{1})-(\alpha_{\ell}({u}_{1})+\beta_{\ell}({u}_{1}))u_{\ell} \, ,
   \ \ell=3,\ldots,8 \, .
\end{eqnarray*}
For the definition of variables $V$, $[Ca]_{i}$, $h$, $j$, $m$,
$d$, $f$, $X$, parameters $I_{st}$, $C_{m}$, $g_{Na}$, $g_{si}$,
$g_{Kp}$, $g_{b}$, $[Na]_{0}$, $[Na]_{i}$, $[K]_{0}$, $[K]_{i}$,
$PR_{NaK}$, $E_{b}$, $T$, constants $c_{1},\ldots,c_{5}$,
functions $g_{K}$, $E_{Na}$, $E_{K}$, $E_{K1}$, $E_{Kp}$,
$K1_{\infty}$, $X_{i}$, $Kp$, $\alpha_{\ell}$, $\beta_{\ell}$,
default values of parameters and initial values of variables in
the Luo-Rudy I model, the reader is referred to
\cite{BichirCL_LuoRudyI}. The reader is also referred to
\cite{BichirCL_LivshitzRudy} for the continuity of the model, and
to \cite{BichirCL_CLB3} for the treatment of the vector field
$\mathcal{F}$ singularities. $\mathcal{F}$ is of class $C^{2}$ on
the domain of interest.

\section{The one - parameter dynamical
system associated to the Luo - Rudy I model}
\label{BichirCL_sect_III}

We performed the study of the dynamical system associated with the
Cauchy problem (\ref{BichirCL_e01}), (\ref{BichirCL_e02}) by
considering only the parameter $\eta_{1}=I_{st}$ and fixing the
rest of parameters. Denote $\lambda=\eta_{1}=I_{st}$ and
$\eta_{\ast}$ the vector of the fixed values of
$\eta_{2},\ldots,\eta_{13}$. Let $F:\mathbb{R} \times M
\rightarrow M$, $F(\lambda,u)=\mathcal{F}(\lambda,\eta_{\ast},u)$,
$F = (F_{1},\ldots,F_{8})$.

Consider the dynamical system associated with the Cauchy problem
(\ref{BichirCL_e01}), (\ref{BichirCL_e06}), where
\begin{equation}
\label{BichirCL_e06}
   \frac{du}{dt}=F(\lambda,u) \, .
\end{equation}

The equilibrium points of this problem are solutions of the
equation
\begin{equation}
\label{BichirCL_e07}
   F(\lambda,u)=0 \, .
\end{equation}
The existence of the solutions and the number were established by
graphical representation in \cite{BichirCL_CLB3}, for the domain
of interest. The equilibrium curve (the bifurcation diagram) was
obtained in \cite{BichirCL_CLB3}, via an arc-length-continuation
method (\cite{BichirCL_CLB10}) and Newton's method
(\cite{BichirCL_Gir_Rav1996}), starting from a solution obtained
by solving a nonlinear least-squares problem
(\cite{BichirCL_CLB10}) for a value of $\lambda$ for which the
system has one solution. In \cite{BichirCL_CLB3}, the results are
obtained by reducing (\ref{BichirCL_e07}) to a system of two
equations in $(u_{1},u_{2})$ $=$ $(V, [Ca]_{i})$. Here, we used
directly (\ref{BichirCL_e07}).

\section{Extended system method for limit cycles}
\label{BichirCL_sect_IV}

The extended system in $(\lambda,T,u)$
\begin{eqnarray}
   & \ & \frac{du}{d \tau}-TF(\lambda,u) = 0 \, ,
         \nonumber \\
   & \ & u(0)-u(1)=0 \, ,
         \label{BichirCL_e23} \\
   & \ & \int\limits_{0}^{1} <u(t), \frac{d w(t)}{d t}> dt = 0 \,
         \nonumber
\end{eqnarray}
is introduced, in \cite{BichirCL_Kuznetsov},
\cite{BichirCL_Doedel}, \cite{BichirCL_MATCONT}, to locate limit
cycles of a general problem (\ref{BichirCL_e01}),
(\ref{BichirCL_e06}). $T$ is the unknown period of the cycle. $w$
is a component of a known reference solution
$(\hat{\lambda},\hat{T},w)$ of (\ref{BichirCL_e23}). The system
(\ref{BichirCL_e23}) becomes determined in a continuation process.

In our case, $<u,v>$ $=$ $\sum\limits_{i=1}^{8}u_{i}v_{i}$ and $\|
u \|$ $=$ $\sqrt{\sum\limits_{i=1}^{8} u_{i}^{2}}$ for $u,v$ $\in$
$\mathbb{R}^{8}$.

In order to approximate and solve (\ref{BichirCL_e23}) by finite
element method time approximation (FEM), let us obtain the weak
form of (\ref{BichirCL_e23}) in the sequel.

Let
\begin{eqnarray*}
   & \ & X=\{ x \in L^{2}(0,1;\mathbb{R}^{8}); \,
         \frac{d x}{d t} \in L^{2}(0,1;\mathbb{R}^{8}),
         \nonumber \\
   & \ & \qquad x=(x_{1},\ldots,x_{8}), \,
         x_{i}(0)=x_{i}(1), \, i=1,\ldots,8  \} \, .
         \label{BichirCL_e39} \\
   & \ & V=\{ v \in L^{2}(0,1;\mathbb{R}); \,
         \frac{d v}{d t} \in L^{2}(0,1;\mathbb{R}), \,
         v(0)=v(1)  \} \, .
         \nonumber
\end{eqnarray*}

The weak form of (\ref{BichirCL_e23}) is the problem in
$(\lambda,T,u)$ $\in$ $\mathbb{R}$ $\times$ $\mathbb{R}$ $\times$
$X$
\begin{eqnarray}
   & \ & \int\limits_{0}^{1} u_{i}(\tau) \ \frac{dv(\tau)}{d \tau} \ d \tau
         + T \int\limits_{0}^{1} F_{i}(\lambda,u(\tau))v(\tau) d \tau = 0 \, ,
         \nonumber \\
   & \ & \qquad \forall v \in V, \, i=1,\ldots,8 \, ,
         \label{BichirCL_e25} \\
   & \ & \int\limits_{0}^{1} <u(t), \frac{d w(t)}{d t}> dt = 0 \, .
         \nonumber
\end{eqnarray}

\section{Arc-length-continuation method for (\ref{BichirCL_e25})}
\label{BichirCL_sect_V}

Following the usual practice (\cite{BichirCL_Keller1},
\cite{BichirCL_Keller2}, \cite{BichirCL_MATCONT},
\cite{BichirCL_Doedel}, \cite{BichirCL_Gir_Rav1996},
\cite{BichirCL_CLB10}, \cite{BichirCL_Govaerts},
\cite{BichirCL_Cl_MatContM}, \cite{BichirCL_Kuznetsov},
\cite{BichirCL_ParkerChua}, \cite{BichirCL_Sey1},
\cite{BichirCL_Sey}, \cite{BichirCL_Sey2}), we also use an
arc-length-continuation method in order to formulate an algorithm
to solve (\ref{BichirCL_e25}) approximatively.

Glowinski (\cite{BichirCL_CLB10}, following H.B.Keller
\cite{BichirCL_Keller1}, \cite{BichirCL_Keller2}) and Doedel
(\cite{BichirCL_Doedel}, where also Keller's name is cited) chose
a continuation equation written in our case as
\begin{equation}
\label{BichirCL_e26}
   \int\limits_{0}^{1} \| \frac{d u (t)}{d s} \|^{2} d t
      + (\frac{d T }{d s})^{2}
      + (\frac{d \lambda }{d s})^{2} = 1 \, ,
\end{equation}
where $s$ is the curvilinear abscissa.

Let $(\lambda^{0},u^{0})$ be a Hopf bifurcation point, $\pm
\beta^{0} i$ a pair of purely imaginary eigenvalues of of the
Jacobian matrix $D_{u}F(\lambda^{0},u^{0})$, and a nonzero complex
vector $g^{0}=g_{r}^{0}+ig_{i}^{0}$. $(\lambda^{0},u^{0})$ is
located on the equilibrium curve during a continuation procedure
using some test functions (\cite{BichirCL_Kuznetsov},
\cite{BichirCL_Govaerts}, \cite{BichirCL_MATCONT}).
$(\lambda^{0},\beta^{0},u^{0},g_{r}^{0},g_{i}^{0})$ $\in$
$\mathbb{R}$ $\times$ $\mathbb{R}$ $\times$ $\mathbb{R}^{8}$
$\times$ $\mathbb{R}^{8} \times \mathbb{R}^{8}$ is the solution of
the extended system (\cite{BichirCL_Sey1}, \cite{BichirCL_Sey},
\cite{BichirCL_Sey2})
\begin{equation}
\label{BichirCL_e24}
   \left[\begin{array}{l}
      F(\lambda,u) \\
      D_{u}F(\lambda,u)g_{r}+\beta g_{i}\\
      D_{u}F(\lambda,u)g_{i}-\beta g_{r}\\
      g_{r,k} - 1\\
      g_{i,k}
      \end{array}\right] = 0 \, ,
\end{equation}
where $k$ is a fixed index of $g_{r}$ and of $g_{i}$, $1 \leq k
\leq 8$.

To solve (\ref{BichirCL_e25}), the extended system formed by
(\ref{BichirCL_e25}) and (\ref{BichirCL_e26}), parametrized by
$s$, was considered. Let $\triangle s$ be an arc-length step and
$\lambda^{n} \cong u(\lambda \triangle s)$, $T^{ \, n} \cong T(n
\triangle s)$, $u^{n} \cong u(n \triangle s)$. We have the
algorithm (following the cases from \cite{BichirCL_CLB10},
\cite{BichirCL_Doedel}, \cite{BichirCL_Sey},
\cite{BichirCL_Sey2}):

1. take the Hopf bifurcation point $(\lambda^{0},u^{0})$ and $T^{
\, 0}=2 \pi / \beta^{0}$; retain $g_{r}^{0}$, $g_{i}^{0}$;

2. for $n = 0$, $(\lambda^{1},T^{ \, 1},u^{1})$ $\in$ $\mathbb{R}$
$\times$ $\mathbb{R}$ $\times$ $X$ is obtained
(\cite{BichirCL_Doedel}, \cite{BichirCL_Sey2}) by
(\ref{BichirCL_e27}),
\begin{eqnarray}
   & \ & \int\limits_{0}^{1} \sum\limits_{i=1}^{8}
         u^{1}_{i}(t) \: \frac{d \phi_{i}(t)}{d t} \ dt = 0 \, \, ,
         \label{BichirCL_e30}
\end{eqnarray}
and
\begin{eqnarray}
   & \ & \int\limits_{0}^{1} \sum\limits_{i=1}^{8}
         (u^{1}_{i}(t)-u^{0}_{i}(t))
        \phi_{i}(t) \ d t
        = \triangle s \, ,
         \label{BichirCL_e32}
\end{eqnarray}
where
\begin{equation}
\label{BichirCL_e29}
   \phi(t)=\sin (2 \pi t) g_{r}^{0}
      + \cos (2 \pi t) g_{i}^{0} \, ,
\end{equation}
using Newton's method with the initial iteration
\begin{equation}
\label{BichirCL_e33}
   (u^{1})^{0}(t)=u^{0}+\triangle s \, \phi(t) \, ,
      \quad (T^{ \, 1})^{0}=T^{ \, 0} \, ,
      \quad (\lambda^{1})^{0}=\lambda^{0} \, .
\end{equation}

 3. for $n \geq 1$, assuming that $(\lambda^{n-1}$, $T^{
\, n-1}$, $u^{n-1})$, $(\lambda^{n}$, $T^{ \, n}$, $u^{n})$ are
known, $(\lambda^{n+1}$, $T^{ \, n+1}$, $u^{n+1})$ $\in$
$\mathbb{R}$ $\times$ $\mathbb{R}$ $\times$ $X$ is obtained by
(\ref{BichirCL_e27}), (\ref{BichirCL_e31}), and
(\ref{BichirCL_e28}), where
\begin{eqnarray}
   & \ & \int\limits_{0}^{1}
         u_{i}^{n+1}(\tau) \ \frac{d v(\tau)}{d \tau} \ d \tau
         + T^{ \, n+1} \int\limits_{0}^{1}
         F_{i}(\lambda^{n+1},u^{n+1}(\tau))v(\tau) d \tau = 0 \, ,
         \label{BichirCL_e27} \\
   & \ & \qquad \forall v \in V, \, i=1,\ldots,8 \, ,
         \nonumber
\end{eqnarray}

\begin{eqnarray}
   & \ & \int\limits_{0}^{1} \sum\limits_{i=1}^{8}
         u^{n+1}_{i}(t) \: \frac{d u^{n}_{i}(t)}{d t} \ dt = 0 \, \, ,
         \label{BichirCL_e31}
\end{eqnarray}

\begin{eqnarray}
   & \ & \int\limits_{0}^{1} \sum\limits_{i=1}^{8}
         (u^{n+1}_{i}(t)-u^{n}_{i}(t))
         \frac{u^{n}_{i}(t)-u^{n-1}_{i}(t)}{\triangle s} \ d t
         \label{BichirCL_e28} \\
   & \ & \qquad + (T^{ \, n+1}-T^{ \, n})\frac{T^{ \, n}-T^{ \, n-1}}{\triangle s}
         + (\lambda^{n+1}-\lambda^{n})\frac{\lambda^{n}-\lambda^{n-1}}{\triangle s}
         = \triangle s \, ,
         \nonumber
\end{eqnarray}
using Newton's method with the initial iteration
\begin{equation}
\label{BichirCL_e40}
   ((\lambda^{n+1})^{0},(T^{ \, n+1})^{0},(u^{n+1})^{0})
   = (\lambda^{n},T^{ \, n},u^{n}). \, .
\end{equation}

\section{Newton's method for the steps of the algorithm from
the end of section \ref{BichirCL_sect_V}} \label{BichirCL_sect_VI}

In (\ref{BichirCL_e28}) ($n \geq 1$), let us denote
$\lambda^{\ast}$ $=$ $\lambda^{n}$, $T^{ \, \ast}$ $=$ $T^{ \,
n}$, $u^{\ast}$ $=$ $u^{n}$, $\lambda^{\ast\ast}$ $=$
$\frac{\lambda^{n}-\lambda^{n-1}}{\triangle s}$, $T^{ \,
\ast\ast}$ $=$ $\frac{T^{ \, n}-T^{ \, n-1}}{\triangle s}$,
$u^{\ast\ast}$ $=$ $\frac{u^{n}-u^{n-1}}{\triangle s}$. We write
(\ref{BichirCL_e27}), (\ref{BichirCL_e30}), (\ref{BichirCL_e32})
(the iteration $n = 0$) in the same general form as
(\ref{BichirCL_e27}), (\ref{BichirCL_e31}), (\ref{BichirCL_e28}).
So denote $u^{\ast}$ $=$ $u^{0}$, $u^{\ast\ast}$ $=$ $\phi$ and
consider $\lambda^{\ast}$ $=$ $\lambda^{0}$, $T^{ \, \ast}$ $=$
$T^{ \, 0}$, $\lambda^{\ast\ast}=0$, $T^{ \, \ast\ast}=0$ in
(\ref{BichirCL_e28}) and consider $u^{\ast}$ $=$ $u^{0}$ $=$
$\phi$ in (\ref{BichirCL_e31}).

Each step of the algorithm at the end of section
\ref{BichirCL_sect_V}, given $(\lambda^{\ast}$, $T^{ \, \ast}$,
$u^{\ast})$, $(\lambda^{\ast\ast}$, $T^{ \, \ast\ast}$,
$u^{\ast\ast})$, calculates $(\lambda^{n+1}$, $T^{ \, n+1}$,
$u^{n+1})$ $\in$ $\mathbb{R}$ $\times$ $\mathbb{R}$ $\times$ $X$,
$n \geq 0$, by (\ref{BichirCL_e27}),
\begin{eqnarray}
   & \ & \int\limits_{0}^{1} \sum\limits_{i=1}^{8}
         u^{n+1}_{i}(t) \, \frac{d u^{\ast}_{i}(t)}{d t} \ dt = 0 \, \, ,
         \label{BichirCL_e34}
\end{eqnarray}
and
\begin{eqnarray}
   & \ & \int\limits_{0}^{1} \sum\limits_{i=1}^{8}
         (u^{n+1}_{i}(t)-u^{\ast}_{i}(t))u^{\ast\ast}_{i}(t) \ d t
         \label{BichirCL_e35} \\
   & \ & \qquad + (T^{ \, n+1}-T^{ \, \ast})T^{ \, \ast\ast}
         + (\lambda^{n+1}-\lambda^{\ast})\lambda^{\ast\ast}
         = \triangle s \, .
         \nonumber
\end{eqnarray}

Newton's method applied (\ref{BichirCL_e27}), (\ref{BichirCL_e34})
and (\ref{BichirCL_e35}), for $n \geq 0$, leads to: let
$((\lambda^{1})^{0}$, $(T^{ \, 1})^{0}$, $(u^{1})^{0})$, given by
(\ref{BichirCL_e33}), be an initial iteration ($m = 0$) if $n =
0$; let $((\lambda^{n+1})^{0}$, $(T^{ \, n+1})^{0}$,
$(u^{n+1})^{0})$, given by (\ref{BichirCL_e40}), be an initial
iteration ($m = 0$) if $n \geq 1$; calculate $(\lambda^{n+1}$,
$T^{ \, n+1}$, $u^{n+1})$ as the solution of the algorithm: for $m
\geq 0$, $((\lambda^{n+1})^{m+1}$, $(T^{ \, n+1})^{m+1}$,
$(u^{n+1})^{m+1})$ $=$ $(\lambda^{m+1}$, $T^{ \, m+1}$, $u^{m+1})$
$\in$ $\mathbb{R}$ $\times$ $\mathbb{R}$ $\times$ $X$ is obtained
by
\begin{eqnarray}
   & \ & \int\limits_{0}^{1} u_{i}^{m+1}(\tau) \ \frac{dv(\tau)}{d \tau} \ d \tau
         + T^{ \, m+1} \int\limits_{0}^{1} F_{i}(\lambda^{m},u^{m}(\tau))v(\tau) d \tau  \, ,
         \nonumber \\
   & \ & + T^{ \, m} \int\limits_{0}^{1} DF_{i}(\lambda^{m},u^{m}(\tau))
         (\lambda^{m+1},u^{m+1}(\tau))v(\tau) d \tau
         \label{BichirCL_e36} \\
   & \ & = T^{ \, m} \int\limits_{0}^{1} DF_{i}(\lambda^{m},u^{m}(\tau))
         (\lambda^{m},u^{m}(\tau))v(\tau) d \tau \, ,
         \nonumber \\
   & \ & \qquad \forall v \in V, \, i=1,\ldots,8 \, ,
         \nonumber
\end{eqnarray}
\begin{eqnarray}
   & \ & \int\limits_{0}^{1} \sum\limits_{i=1}^{8}
         u^{m+1}_{i}(t) \, \frac{d u^{\ast}_{i}(t)}{d t} \ dt = 0 \, \, ,
         \label{BichirCL_e37}
\end{eqnarray}
\begin{eqnarray}
   & \ & \int\limits_{0}^{1} \sum\limits_{i=1}^{8}
         u^{m+1}_{i}(t)u^{\ast\ast}_{i}(t) \ d t
         + T^{ \, m+1}T^{ \, \ast\ast}
         + \lambda^{m+1}\lambda^{\ast\ast}
         \label{BichirCL_e38} \\
   & \ & \qquad = \int\limits_{0}^{1} \sum\limits_{i=1}^{8}
         u^{\ast}_{i}(t)u^{\ast\ast}_{i}(t) \ d t
         + T^{ \, \ast}T^{ \, \ast\ast}
         + \lambda^{\ast}\lambda^{\ast\ast}
         + \triangle s \, .
         \nonumber
\end{eqnarray}

\section{Approximation of problem (\ref{BichirCL_e36}),
(\ref{BichirCL_e37}), (\ref{BichirCL_e38}) by finite element
method time approximation} \label{BichirCL_sect_VII}

In order to perform this approximation, let us divide the interval
$[0,1]$ in $N+1$ subintervals $K=K_{j}=[t_{j},t_{j+1}]$, $0 \leq j
\leq N$, where $0 = t_{0} < t_{1} < \ldots < t_{N+1} = 1$. The
sets $K$ represent a triangulation $\mathcal{T}_{h}$ of $[0,1]$.

Let us approximate the spaces $V$ and $X$ by the spaces
\begin{eqnarray*}
   & \ & V_{h}=\{ v : [0,1] \rightarrow \mathbb{R}; \,
         v \in C[0,1], \,
         v(0)=v(1), \,
         v|_{K} \in P_{k}(K), \, \forall K \in \mathcal{T}_{h}
         \} \, ,
         \label{BichirCL_e41} \\
   & \ & X_{h}=\{ x  : [0,1] \rightarrow \mathbb{R}^{8}; \,
         x=(x_{1},\ldots,x_{8}), \,
         x_{i} \in V_{h}, \, i=1,\ldots,8  \} \, ,
         \nonumber
\end{eqnarray*}
respectively, where $P_{k}(K)$ is the space of polynomials in $t$
of degree less than or equal to $k$ defined on $K$, $k \geq 2$.

Let $k = 2$. An element $K \in \mathcal{T}_{h}$ has three nodal
points. To obtain a function $u_{h}$ $\in$ $X_{h}$ reduces to
obtain a function $v_{h}$ $\in$ $V_{h}$. In order to obtain a
function $v_{h}$ $\in$ $V_{h}$, we use a basis of functions of
$V_{h}$. Let $J_{K}=\{ 1,2,3 \}$ be the local numeration for the
nodes of $K$, where $1, 3$ correspond to $t_{j}, t_{j+1}$
respectively and $2$ corresponds to a node between $t_{j}$ and
$t_{j+1}$. Let $\{ \psi_{i}, i \in J_{K} \}$ be the local
quadratic basis of functions on $K$ corresponding to the local
nodes. Let $J=\{ 1, \ldots, 2N+1 \}$ be the global numeration for
the nodes of $[0,1]$. The two numerations are related by a matrix
$L$ whose elements are the elements $j \in J$. Its rows are
indexed by the elements $K \in \ \mathcal{T}_{h}$ (by the number
of the element $K$ in a certain fixed numeration with elements
from the set $\{ 1,\ldots,N \}$) and its columns, by the local
numeration $i \in J_{K}$, that is $j = L(K,i)$. A function $v_{h}$
$\in$ $V_{h}$ is defined by its values $v_{j}$ from the nodes $j
\in J$,
\begin{equation}
\label{BichirCL_e42}
   v_{h}(t)
      = \sum_{K \in \ \mathcal{T}_{h}}  \;
     \sum_{i \in \ J_{K}, \; j = L(K,i)}
        v_{j} \, \psi_{i}(t) \, ,
\end{equation}
and a function $u_{h}$ $\in$ $X_{h}$ is defined by its values
$u_{j}$ from the nodes $j \in J$,
\begin{equation}
\label{BichirCL_e43}
   u_{h}(t)
      = \sum_{K \in \ \mathcal{T}_{h}}  \;
     \sum_{i \in \ J_{K}, \; j = L(K,i)}
        u_{j} \, \psi_{i}(t) \, .
\end{equation}
So, an unknown function $u_{h}$ $=$ $((u_{h})_{1}$, $ \ldots$,
$(u_{h})_{8})$ is reduced to the unknowns $u_{j}$, $u_{j}$ $=$
$((u_{j})_{1}$, $\ldots$, $(u_{j})_{8})$, $j \in J$.

In (\ref{BichirCL_e36}), (\ref{BichirCL_e37}),
(\ref{BichirCL_e38}), approximate $(\lambda^{m+1}$, $T^{ \, m+1}$,
$u^{m+1})$ $\in$ $\mathbb{R}$ $\times$ $\mathbb{R}$ $\times$ $X$
by $(\lambda^{m+1}$, $T^{ \, m+1}$, $u^{m+1}_{h})$ $\in$
$\mathbb{R}$ $\times$ $\mathbb{R}$ $\times$ $X_{h}$. Taking
$u^{m+1}_{h} = u_{h}$, $u_{h}$ given by (\ref{BichirCL_e43}), and
$v = \psi_{\ell}$, for all $\ell \in J_{K}$, for all $K \in \
\mathcal{T}_{h}$, we obtain the discrete variant of problem
(\ref{BichirCL_e36}), (\ref{BichirCL_e37}), (\ref{BichirCL_e38})
as the following problem in $(\lambda$, $T$, $u_{1}$, $\ldots$,
$u_{2N+1})$ $\in$ $\mathbb{R}$ $\times$ $\mathbb{R}$ $\times$
$\mathbb{R}^{8 \cdot (2N+1)}$, written suitable for the assembly
process,
\begin{eqnarray}
   & \ & \sum_{K \in \ \mathcal{T}_{h}} \{ \;
         \sum_{i \in \ J_{K}, \; j = L(K,i)}
         (u_{j}^{m+1})_{n} \int\limits_{K}
         \psi_{i}(\tau) \: \frac{\psi_{\ell} \, (\tau)}{d \tau} \; d \tau
         \nonumber \\
   & \ & \qquad + \, T^{ \, m+1} \int\limits_{K}
         F_{i}(\lambda^{m},u^{m}_{h}(\tau)) \; \psi_{\ell} \; d \tau
         \nonumber \\
   & \ & \qquad + \, T^{ \, m} \sum_{i \in \ J_{K}, \; j = L(K,i)}
         <u_{j}^{m+1}, \int\limits_{K}
         (D_{u}F_{n}(\lambda^{m},u^{m}_{h}(\tau))
         \psi_{i}(\tau)) \, \psi_{\ell} \; d \tau >
         \nonumber \\
   & \ & \qquad + \lambda^{m+1} \int\limits_{K}
         (D_{\lambda}F_{n}(\lambda^{m},u^{m}_{h}(\tau)))
         \; \psi_{\ell} \; d \tau \; \}
         \label{BichirCL_e44} \\
   & \ & = \, T^{ \, m} \sum_{K \in \ \mathcal{T}_{h}} \{ \;
         \int\limits_{K}
         (D_{u}F_{n}(\lambda^{m},u^{m}_{h}(\tau))
         u^{m}_{h}(\tau)) \, \psi_{\ell} \; d \tau
         \nonumber \\
   & \ & \qquad + \int\limits_{K}
         (D_{\lambda}F_{n}(\lambda^{m},u^{m}_{h}(\tau))
         \lambda^{m}) \; \psi_{\ell} \; d \tau \; \} \, ,
         \, n=1,\ldots,8 \, ,
         \nonumber \\
   & \ & \qquad u_{0}^{m+1} = u_{2N+1}^{m+1} \, ,
         \label{BichirCL_e44bis}
\end{eqnarray}
\begin{eqnarray}
   & \ & \sum_{K \in \ \mathcal{T}_{h}}  \;
         \sum_{i \in \ J_{K}, \; j = L(K,i)} \;
         <
         u_{j}^{m+1}, \int\limits_{K}
         \psi_{i}(\tau) \, \frac{d u^{\ast}_{h}(\tau)}{d \tau}
         \ d \tau > \, = \, 0 \, \, ,
         \label{BichirCL_e45}
\end{eqnarray}
\begin{eqnarray}
   & \ & \sum_{K \in \ \mathcal{T}_{h}}  \;
         \sum_{i \in \ J_{K}, \; j = L(K,i)} \;
         <
         u_{j}^{m+1}, \int\limits_{K}
         \psi_{i}(\tau) \, u^{\ast\ast}_{h}(\tau) \ d \tau >
         \label{BichirCL_e46} \\
   & \ & \qquad + T^{ \, m+1}T^{ \, \ast\ast}
         + \lambda^{m+1}\lambda^{\ast\ast}
        \nonumber \\
   & \ & \qquad = \sum_{K \in \ \mathcal{T}_{h}}
         \int\limits_{K}
         < u^{\ast}_{h}(\tau),  u^{\ast\ast}_{h}(\tau) > \; d \tau
        + T^{ \, \ast}T^{ \, \ast\ast}
         + \lambda^{\ast}\lambda^{\ast\ast}
         + \triangle s \, ,
        \nonumber
\end{eqnarray}
for all $\ell \in J_{K}$, for all $K \in \ \mathcal{T}_{h}$.

\begin{figure}
   \raggedleft
   \includegraphics[height=9cm,width=14cm]{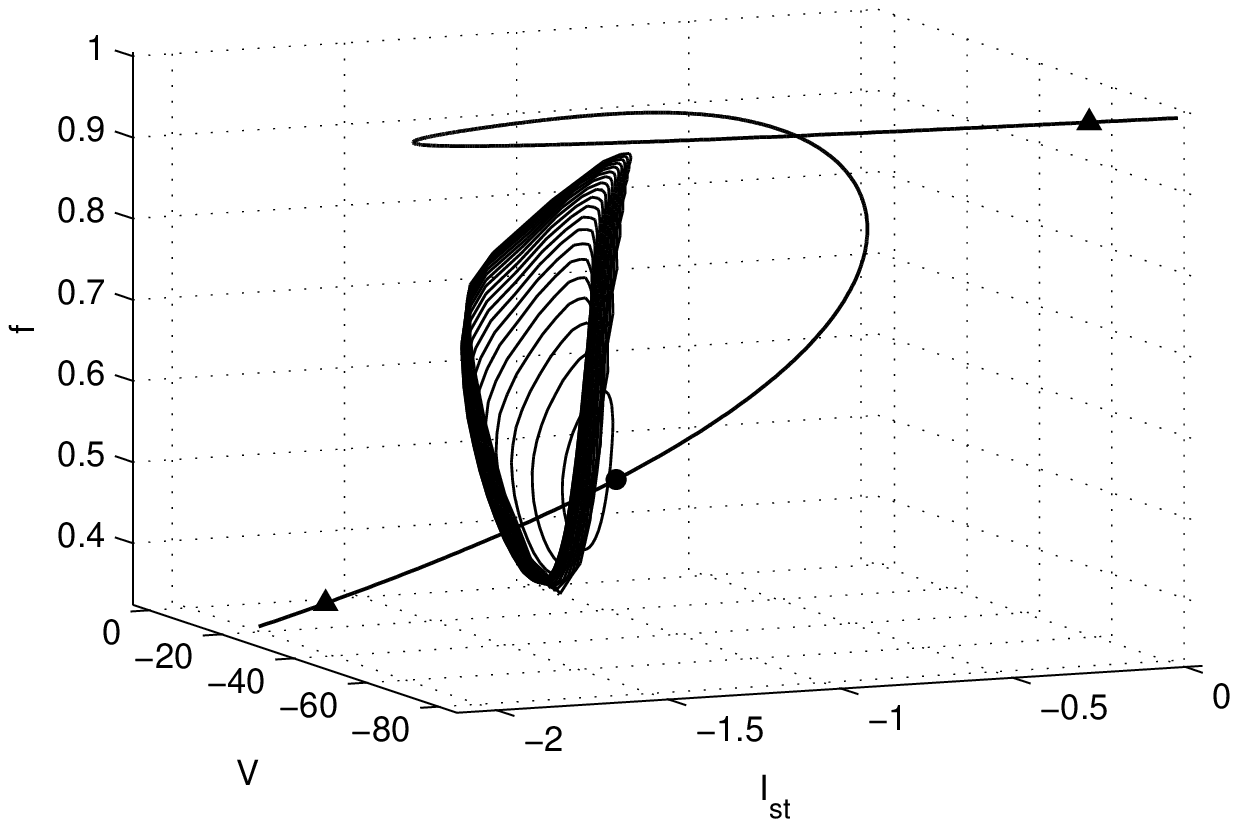}
   \raggedleft
   \includegraphics[height=9cm,width=14cm]{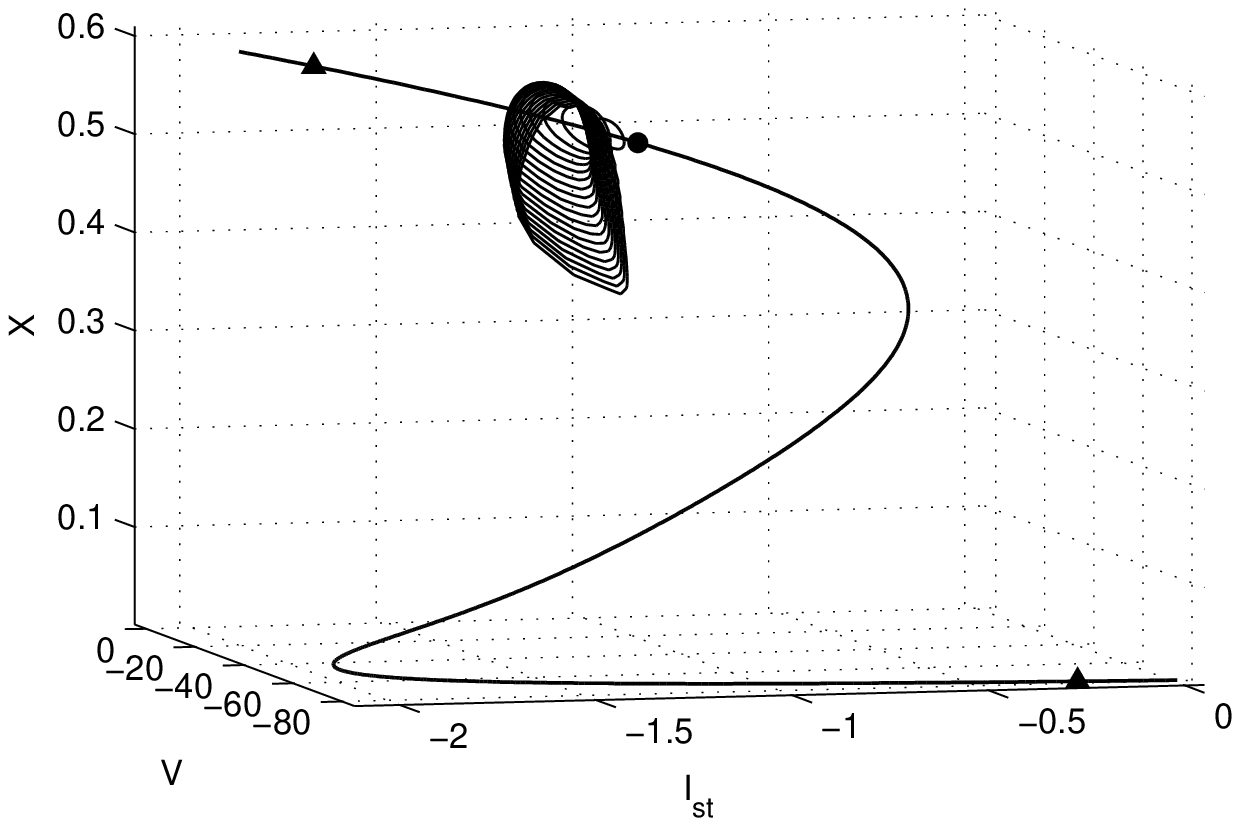}
   \caption{Two projections of limit cycles and of a part of the equilibrium
   curve (marked by "$\blacktriangle$"). The Hopf bifurcation
   point is marked by "$\bullet$".}
   \label{BichirCL_proj_LC}
\end{figure}

\begin{figure}
   \raggedleft
   \includegraphics[height=4cm,width=6cm]{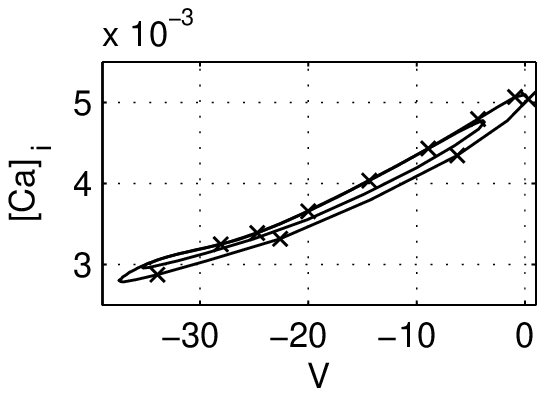}
   \centering
   \includegraphics[height=4cm,width=6cm]{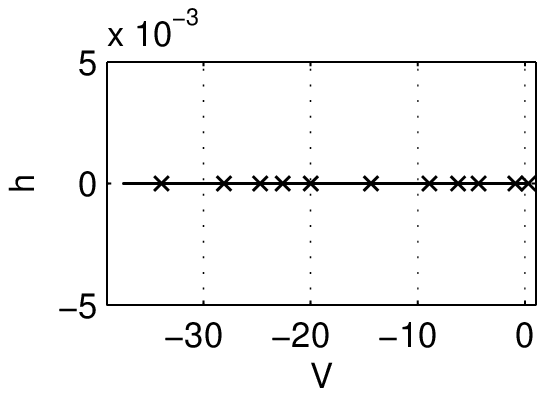}
   \raggedleft
   \includegraphics[height=4cm,width=6cm]{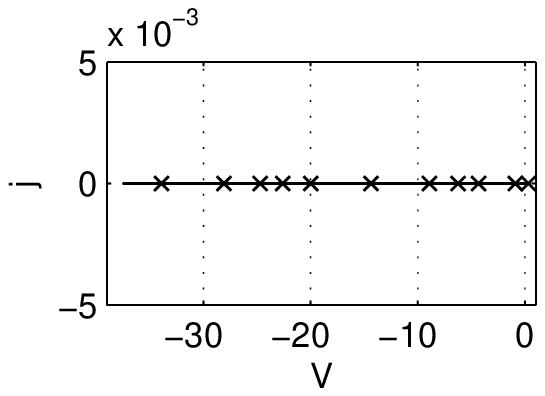}
   \centering
   \includegraphics[height=4cm,width=6cm]{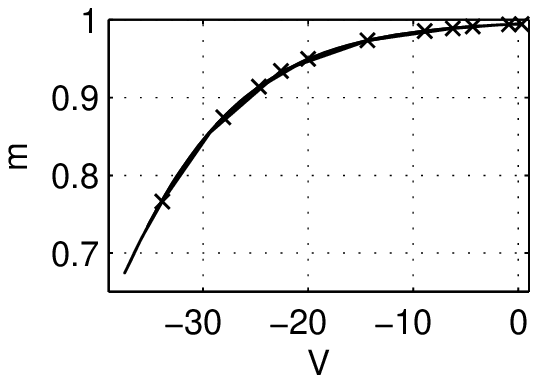}
   \raggedleft
   \includegraphics[height=4cm,width=6cm]{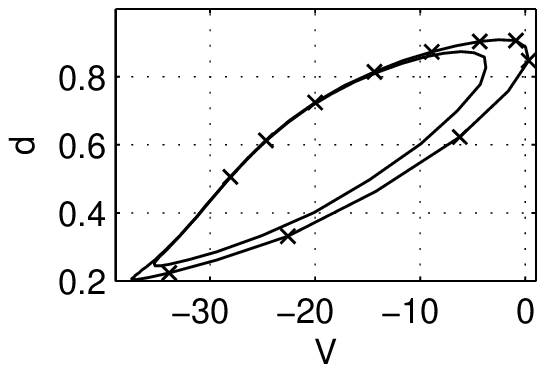}
   \centering
   \includegraphics[height=4cm,width=6cm]{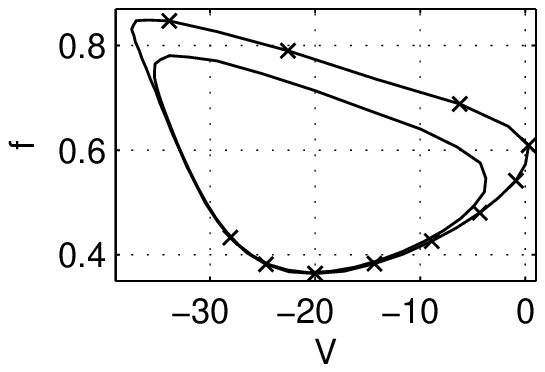}
   \centering
   \includegraphics[height=4cm,width=6cm]{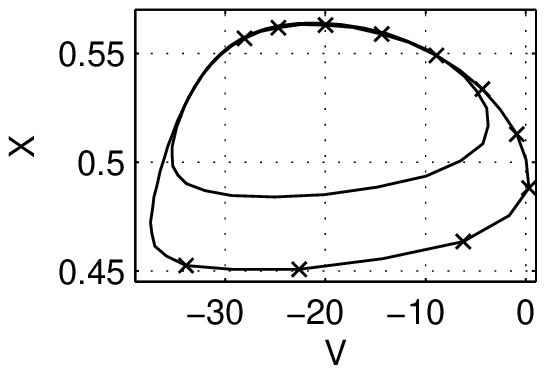}
   \caption{Projections of two limit cycles calculated for
   $I_{st}$ $=$ $-1.2000465026$ and for $I_{st}$ $=$ $-1.2000183729$
   (marked by "x")
   (20 elements, 41 nodes).}
   \label{BichirCL_proj_LC_20elem}
\end{figure}

\section{Numerical results}
\label{BichirCL_sect_X} Based on \cite{BichirCL_Ta_Hug} and on the
computer programs for \cite{BichirCL_CLB1} and
\cite{BichirCL_CLB2}, relations (\ref{BichirCL_e44}),
(\ref{BichirCL_e44bis}), (\ref{BichirCL_e45}),
(\ref{BichirCL_e46}) and the algorithm at the end of section
\ref{BichirCL_sect_V} furnished the numerical results of this
section.

Let $(\lambda^{0},u^{0})$ be the Hopf bifurcation point located
during the construction of the equilibrium curve by a continuation
procedure in \cite{BichirCL_CLB4}. The solution
$(\lambda^{0},\beta^{0},u^{0},g_{r}^{0},g_{i}^{0})$ of
(\ref{BichirCL_e24}), calculated in \cite{BichirCL_CLB4}, is
$\lambda^{0}$ $=$ $-1.0140472901$, $\beta^{0}$ $=$ $0.0162886062$,
$u^{0}$ $=$ $(-24.3132508542$, $0.0034641214$, $0.0$, $0.0$,
$0.9176777444$, $0.5025242162$, $0.4920204612$, $0.5071561613)$,
$g_{r}^{0}$ $=$ $(1.0$, $0.0000468233$, $0.0$, $0.0$,
$0.0093195354$, $0.0198748652$, $-0.0072420216$, $0.0001706577)$,
$g_{i}^{0}$ $=$ $(0.0$, $0.0000029062$, $0.0$, $0.0$,
$-0.0000171311$, $-0.0118192370$, $0.0136415789$,
$-0.0017802907)$. (The eigenvalues of the Jacobian matrix
$D_{u}F(\lambda^{0},u^{0})$, calculated by the QR algorithm, are
$\pm \, 0.0162886062 \, i$, $-8.8611865338$, $-0.1026761869$,
$-0.0647560667$, $-0.0024565181$, $-1.7398266947$,
$-0.2049715178$). These data are considered in the step 1 of the
algorithm at the end of section \ref{BichirCL_sect_V}.

We took $C_{m}=1$, $g_{Na}=23$, $g_{si}=0.09$, $g_{K}=0.282$,
$g_{K1}=0.6047$, $g_{Kp}=0.0183$, $G_{b}=0.03921$, $[Na]_{0}=140$,
$[Na]_{i}=18$, $[K]_{0}=5.4$, $[K]_{i}=145$, $PR_{NaK}=0.01833$,
$E_{b}=-59.87$, $T=310$.

In order to solve (\ref{BichirCL_e25}) numerically by the
algorithm at the end of section \ref{BichirCL_sect_V} and by
(\ref{BichirCL_e44}), (\ref{BichirCL_e44bis}),
(\ref{BichirCL_e45}), (\ref{BichirCL_e46}), we performed
calculations using $\triangle s$ $=$ $1.0$ and 500 iterations in
the continuation process. Integrals $\int\limits_{K} f(\tau) \; d
\tau$ were calculated using Gauss integration formula with three
integration points.

Figure \ref{BichirCL_proj_LC} and \ref{BichirCL_proj_LC_20elem}
present some results obtained using 20 elements $K$ (41 nodes) ($N
= 20$, $J=\{ 1, \ldots, 41 \}$ in section
\ref{BichirCL_sect_VII}). The curves of the projections of the
limit cycles, on the planes indicates in figure, are plots
generated from values calculated in the nodes, corresponding to a
fixed value of the parameter.

Two projections of some limit cycles and of a part of the
equilibrium curve (marked by "$\blacktriangle$") are presented in
figure \ref{BichirCL_proj_LC}. The Hopf bifurcation point is
marked by "$\bullet$".

In figure \ref{BichirCL_proj_LC_20elem}, there are represented the
projections of the plots of two limit cycles calculated for
$I_{st}$ $=$ $-1.2000465026$ (iteration 148) and for $I_{st}$ $=$
$-1.2000183729$ (iteration 248, marked by "x" in figure).

The results obtained are relevant from a biological point of view,
pointing to unstable electrical behavior of the modeled system in
certain conditions, translated into oscillatory regimes such as
early afterdepolarizations
(\cite{BichirCL_Tran_Sato_Yochelis2009}) or self-sustained
oscillations (\cite{BichirCL_CLB3}), which may in turn
synchronize, resulting in life-threatening arrhythmias: premature
ventricular complexes or torsades-de-pointes, degenerating in
rapid polymorphic ventricular tacycardia or fibrillation
(\cite{BichirCL_Sato_Xie_Sovari_sa}).

\textit{Acknowlegdements}: This research was partially supported
from grant PNCDI2 61-010 to M-LF by the Romanian Ministry of
Education, Research, and Innovation.


\begin{thebibliography}{1}
\bibitem{BichirCL_BeelerReuter} G. W. Beeler , H. Reuter, \textit{Reconstruction of the action
potential of ventricular myocardial fibres}, J. Physiol.
\textbf{268}(1977), 177-210.

\bibitem{BichirCL_CLB1} C. L. Bichir, A.Georgescu, \textit{Approximation of pressure
perturbations by FEM}, Scientific Bulletin of the Pite\c{s}ti
University, the Mathematics-Informatics Series, \textbf{9} (2003),
31-36.

\bibitem{BichirCL_CLB2} C. L. Bichir, \textit{A numerical study by FEM and FVM of a problem which
presents a simple limit point}, ROMAI J., \textbf{4}, 2(2008),
45-56, http://www.romai.ro, http://rj.romai.ro.

\bibitem{BichirCL_CLB3} C. L. Bichir, B. Amuzescu, A. Georgescu, M. Popescu,
Ghe. Nistor, I. Svab, M. L. Flonta, A. D. Corlan,
\textit{Stability and self-sustained oscillations in a ventricular
cardiomyocyte model}, submitted to Interdisciplinary Sciences -
Computational Life Sciences, Springer.

\bibitem{BichirCL_CLB4} C. L. Bichir, A. Georgescu, B. Amuzescu, Ghe. Nistor,
M. Popescu, M. L. Flonta, A. D. Corlan, I. Svab, \textit{Limit
points and Hopf bifurcation points for a one - parameter dynamical
system associated to the Luo - Rudy I model}, submitted to
Mathematics and its Applications,
http://www.mathematics-and-its-applications.com.

\bibitem{BichirCL_Cu_Se_St} C. Cuvelier, A.Segal, A.A.van Steenhoven, \textit{Finite
Element Methods and Navier-Stokes Equations}, Reidel, Amsterdam,
1986.

\bibitem{BichirCL_MATCONT} A. Dhooge, W. Govaerts, Yu.A. Kuznetsov, W. Mestrom,
A.M. Riet, B. Sautois, \textit{MATCONT and CL$_{-}$MATCONT:
Continuation toolboxes in MATLAB}, 2006,
http://www.matcont.ugent.be/manual.pdf

\bibitem{BichirCL_Doedel} E. Doedel, \textit{Lecture Notes on Numerical
Analysis of Nonlinear Equations}, 2007,
http://cmvl.cs.concordia.ca/publications/notes.ps.gz, from the
Home Page of the AUTO Web Site, http://indy.cs.concordia.ca/auto/.

\bibitem{BichirCL_AG_Mo_O} A.Georgescu, M.Moroianu, I.Oprea,
\textit{Bifurcation Theory. Principles and Applications}, Applied
and Industrial Mathematics Series, \textbf{1}, University of
Pite\c{s}ti, 1999.

\bibitem{BichirCL_Gibb} W. J. Gibb , M. B. Wagner, M. D. Lesh,
\textit{Effects of simulated potassium blockade on the dynamics of
triggered cardiac activity}, J. theor. Biol \textbf{168}(1994),
245-257.

\bibitem{BichirCL_Gir_Rav1979} V.Girault, P.-A.Raviart, \textit{Finite Element
Approximation of the Navier-Stokes Equations}, Springer, Berlin,
1979.

\bibitem{BichirCL_Gir_Rav1996} V.Girault, P.-A.Raviart, \textit{Finite Element
Methods for Navier-Stokes Equations.Theory and Algorithms},
Springer, Berlin, 1986.

\bibitem{BichirCL_CLB10} R.Glowinski, \textit{Numerical Methods for Nonlinear Variational
Problems}, Springer, New York, 1984.

\bibitem{BichirCL_Govaerts} W.J.F. Govaerts, \textit{Numerical methods for Bifurcations
of Dynamical Equilibria}, SIAM, Philadelphia, 2000.

\bibitem{BichirCL_Cl_MatContM} W. Govaerts, Yu. A. Kuznetsov R. Khoshsiar Ghaziani,
H.G.E. Meijer, \textit{Cl$_{-}$MatContM: A toolbox for
continuation and bifurcation of cycles of maps}, 2008,
http://www.matcont.ugent.be/doc$_{-}$cl$_{-}$matcontM.pdf

\bibitem{BichirCL_Hodgkin_Huxley} A. L. Hodgkin, A. F. Huxley,
\textit{A quantitative description of membrane current and its
application to conduction and excitation in nerve}, J. Physiol.,
\textbf{117} (1952), 500-544.

\bibitem{BichirCL_Keller1} H.B.Keller, \textit{Numerical Solution
of Bifurcation Eigenvalue Problems}, in \textit{Applications in
Bifurcation Theory}, ed. by P.Rabinowitz, Academic, New York,
1977.

\bibitem{BichirCL_Keller2} H.B.Keller, \textit{Global Homotopies and Newton Methods},
in \textit{Recent Advances in Numerical Methods}, ed. by C. de
Boor, G.H.Golub, Academic, New York, 1978.

\bibitem{BichirCL_Kuznetsov} Yu. A. Kuznetsov, \textit{Elements of Applied Bifurcation
Theory}, Springer, New York, 1998.

\bibitem{BichirCL_LivshitzRudy} L. Livshitz, Y. Rudy, \textit{Uniqueness and stability of
action potential models during rest, pacing, and conduction using
problem - solving environment}, Biophysical J., \textbf{97}
(2009), 1265-1276.

\bibitem{BichirCL_LuoRudyI} C.H. Luo, Y. Rudy, \textit{A model of the ventricular cardiac
action potential. Depolarization, repolarization, and their
interaction}, Circ. Res., \textbf{68} (1991), 1501-1526.

\bibitem{BichirCL_Neher_Sakmann} E. Neher, B. Sakmann,  \textit{Single-channel
currents recorded from membrane of denervated frog muscle fibres},
Nature, \textbf{260} (1976), 799-802.

\bibitem{BichirCL_Noble7} D. Noble,  \textit{Modelling the heart:
insights, failures and progress}, Bioessays, \textbf{24} (2002),
1155-1163.

\bibitem{BichirCL_Noble2007} D. Noble,  \textit{From the Hodgkin-Huxley
axon to the virtual hear}, J. Physiol., \textbf{580} (2007),
15-22. Epub 2006 Oct 2005.

\bibitem{BichirCL_ParkerChua} T.S.Parker, L.O.Chua, \textit{Practical Numerical Algorithms for
Chaotic Systems}, Springer, New York, 1989.

\bibitem{BichirCL_Sato_Xie_Sovari_sa} D. Sato, L. H. Xie,
A. A. Sovari, D. X. Tran, N. Morita, F. Xie, H. Karagueuzian, A.
Garfinkel, J. N. Weiss, Z. Qu , \textit{Synchronization of chaotic
early afterdepolarizations in the genesis of cardiac arrhythmias},
Proc. Natl. Acad. Sci. USA, \textbf{106} (2009), 2983-2988. Epub
2009 Feb 2913.

\bibitem{BichirCL_Sey1} R. Seydel, \textit{Numerical computation of branch points in
nonlinear equations}, Numer. Math., \textbf{33} (1979), 339-352.

\bibitem{BichirCL_Sey} R. Seydel, \textit{Nonlinear Computation}, invited
lecture and paper presented at the Distinguished Plenary Lecture
session on Nonlinear Science in the 21st Century, 4th IEEE
International Workshops on Cellular Neural Networks and
Applications, and Nonlinear Dynamics of Electronic Systems,
Sevilla, June, 26, 1996.

\bibitem{BichirCL_Sey2} R. Seydel, \textit{Practical Bifurcation
and Stability Analysis}, Springer, New York, 2010.

\bibitem{BichirCL_Ta_Hug} C. Taylor, T.G. Hughes, \textit{Finite Element Programming of the
Navier-Stokes Equations}, Pineridge Press, Swansea, U.K., 1981.

\bibitem{BichirCL_Tem3} R.Temam, \textit{Navier-Stokes equations.
Theory and numerical analysis}, North-Holland, Amsterdam, 1979.

\bibitem{BichirCL_Tran_Sato_Yochelis2009} D. X, Tran, D. Sato, A. Yochelis,
J. N. Weiss, A. Garfinkel, Z. Qu, \textit{Bifurcation and chaos in
a model of cardiac early afterdepolarizations}, Phys. Rev. Lett.,
\textbf{102:258103} (2009). Epub 252009 Jun 258125.

\bibitem{BichirCL_Van_der_Pol_van_der_Mark} B. Van der Pol,
J. Van der Mark, \textit{The heartbeat considered as a relaxation
oscillation and an electrical model of the heart}, Phil. Mag.
(suppl.), \textbf{6} (1928), 763-775.

\end{thebibliography}

\end{document}